\documentclass{amsart}
\usepackage{url,amsmath,amssymb, graphicx}
\newcommand{\FriCAS}{{\tt FriCAS}}
\newcommand{\sbcl}{{\tt SBCL}}
\newcommand{\MMA}{{\tt Mathematica}}
\newcommand{\Guess}{{\tt Guess}}
\newcommand{\GuessKauers}{{\tt Guess.m}}
\newcommand{\GuessRat}{{\tt GuessRat}}
\newcommand{\GuessHolo}{{\tt GuessHolo}}
\newcommand{\Rate}{{\tt Rate.m}}
\newcommand{\GFUN}{{\tt GFUN}}
\newcommand{\qGeneratingFunctions}{{\tt qGeneratingFunctions.m}}
\newcommand{\SuperSeeker}{{\tt SuperSeeker}}

\newcommand{\spad}[1]{{\tt #1}}

\newcommand{\OEIS}[1]{{\tt #1}}

\DeclareMathOperator{\ord}{ord}
\DeclareMathOperator{\D}{D}
\DeclareMathOperator{\dom}{dom}
\DeclareMathOperator{\rank}{rank}
\DeclareMathOperator{\Span}{span}
\DeclareMathOperator{\defect}{defect}
\DeclareMathOperator{\critical}{critical}

\newcommand{\Mat}[1]{\ensuremath{\mathbf{#1}}}             
\newcommand{\Dfn}[1]{\emph{#1}}                            

\newtheorem{thm}{Theorem}[section]
\newtheorem{lem}[thm]{Lemma}
\newtheorem{cor}[thm]{Corollary}

\newtheorem*{prbSeries}{Rational Interpolation Problem, Series Variant}
\newtheorem*{prbSequence}{Rational Interpolation Problem, Sequence
  Variant}

\theoremstyle{definition}
\newtheorem{dfn}[thm]{Definition}

\newcommand{\shrink}[2][0pt]
{\hbox to #1{\hss #2\hss}}

\newcommand{\mshrinksub}[2]{\hbox to 0pt{$#1#2$\hss}}

\let\size=\abs

\newcommand{\qbinom}[2]{\genfrac{[}{]}{0pt}{}{#1}{#2}}

\author{Waldemar Hebisch}
\address{Institute of Mathematics, Wroc{\l}aw University}
\email{waldemar.hebisch@math.uni.wroc.pl}
\urladdr{http://www.math.uni.wroc.pl/~hebisch/}

\author{Martin Rubey}

\thanks{Research partially supported by the Austrian Science Foundation FWF,
  grant S9607-N13, in the framework of the National Research Network \lq\lq
  Analytic Combinatorics and Probabilistic Number Theory\rq\rq.}
\address{Institut f\"ur Algebra, Zahlentheorie und Diskrete Mathematik, Leibniz
  Universit\"at Hannover, Welfengarten 1, D-30167 Hannover, Germany}
\email{martin.rubey@math.uni-hannover.de}
\urladdr{http://www.iazd.uni-hannover.de/~rubey/}

\title{Extended Rate, more GFUN}

\usepackage[breaklinks=true,hyperref]{hyperref}
\begin{document}
\maketitle
\begin{abstract}
  We present a software package that guesses formulae for sequences of, for
  example, rational numbers or rational functions, given the first few terms.
  We implement an algorithm due to Bernhard Beckermann and George Labahn,
  together with some enhancements to render our package efficient.  Thus we
  extend and complement Christian Krattenthaler's program \Rate, the parts
  concerned with guessing of Bruno Salvy and Paul Zimmermann's \GFUN, the
  univariate case of Manuel Kauers' \GuessKauers\ and Manuel Kauers' and
  Christoph Koutschan's \qGeneratingFunctions.
\end{abstract}

\section{Introduction}
For some a brain-teaser, for others one step in proving their next theorem:
given the first few terms of a sequence of, say, integers, what is the next
term, what is the general formula?  Of course, no unique solution exists, but,
by Occam's razor, we will prefer a \lq simple\rq\ formula over a more \lq
complicated\rq\ one.  In this article we present a new package that aims at
finding such a simple formula, written for the computer algebra
system \FriCAS.\footnote{\FriCAS\ is freely available at
  \url{http://fricas.sourceforge.net}.} 

Some sequences are very easy to \lq guess\rq, like
\begin{equation}\label{eq:nsquared}
0,1,4,9,\dots,
\end{equation}
or
\begin{equation}\label{eq:Fibonacci}
1,1,2,3,5,\dots.
\end{equation}
Others are a little harder, for example
\begin{equation}\label{eq:sumfact}
0,1,3,9,33,\dots.
\end{equation}
Of course, at times we might want to guess a formula for a sequence of
polynomials or rational functions, too:
\begin{equation}\label{eq:binomial}
1,1+q+q^2,(1+q+q^2)(1+q^2),(1+q^2)(1+q+q^2+q^3+q^4)\dots,
\end{equation}
or
\begin{equation}\label{eq:qExpRat}
  \frac{1-2q}{1-q}, 1-2q,(1-q)(1-2q)^3,(1-q)^2(1-2q)(1-2q-2q^2)^3,\dots
\end{equation}

Fortunately, with the right tool, it is a matter of a moment to figure out
formulae for all of these sequences.  In this article we describe a computer
program that encompasses well known techniques and adds new ideas that we hope
to be effective.  In particular, we generalise both Christian Krattenthaler's
program \Rate~\cite{RATE}, and the guessing functions present in \GFUN\ written
by Bruno Salvy and Paul Zimmermann~\cite{GFUN}, and in \qGeneratingFunctions\
by Manuel Kauers and Christoph Koutschan~\cite{RISC3145, MR2511667}.  With a little manual
aid, we can guess multivariate formulae as well, along the lines of Doron
Zeilberger's programs \GuessRat\ and \GuessHolo~\cite{Zeilberger2005,
  Zeilberger2007}, or Manuel Kauers' program \GuessKauers~\cite{RISC09-07}.
All these programs, as well as the one presented here, try to \emph{compute} a
function that yields the terms when evaluated at $0$, $1$, $2$, and so on.  We
describe this computational approach in more detail beginning in
Section~\ref{sec:background}.

A completely different idea is pursued by \emph{The online encyclopedia of
  integer sequences} of Neil Sloane~\cite{OEIS}.  There, you can enter a
sequence of integers and chances are good that the website will respond with
one or more likely matches.  However, the approach taken is quite different
from ours: the encyclopedia keeps a list of currently roughly $160,000$
sequences, entered more or less manually, and it compares the given sequence
with each one of those.  Besides that, there is an email service called
\SuperSeeker\ that tries some transformations on the given sequence to find a
match in the database.  Furthermore it tries some programs in the spirit of
\Rate\ and \GFUN\ to find a formula, although with a time limit, i.e., it gives
up when too much time has elapsed.

Thus, the two approaches complement each other: For example, there are
sequences where no simple formula is likely to exist, and which can thus be
found only in the encyclopedia.  On the other hand, there are many sequences
that have not yet found their way into the encyclopedia, but can be guessed in
easily by your computer.

In Section~\ref{sec:function-classes} we outline the capabilities of our
package.  In Section~\ref{sec:options} we describe the
most important options that modify the behaviour of the functions.
A very brief description of the
algorithms used and the efficiency problems encountered is given in
Section~\ref{sec:background} and thereafter.

\section{History}

On the historical side, we remark that already in 1964, Malcolm Pivar and Mark
Finkelstein~\cite{PivarFinkelstein1964} implemented a program to identify
sequences given their first few terms, see also Paul W.
Abrahams~\cite{Abrahams}. One interesting feature of their program was the
ability to deal with exceptions to a rule: their program would apply, for
example, the differencing operator, until most of the terms would be equal. In
a second step, it would then locate the exceptions to the rule and try to guess
formulae for the positions and for the values of the exceptions.

The first edition \cite{Sloane1973} of \lq A Handbook of Integer Sequences\rq\
by Neil Sloane appeared in 1973.  In 1992, Fran{\c{c}}ois Bergeron and Simon
Plouffe~\cite{BergeronPlouffe1992} explored the idea of applying various
transformations to the given sequence, for example series reversion.  They then
used Pad\'e approximation to see whether the result might be rational.  In the
same article an experimental program to check for \lq constructible
differentially finite\rq\ series is also briefly described, but it seems that
it was not so successful.

In Physics, Richard Brak, G. S. Joyce, Michael E. Fisher, Helen Au-Yang,
Anthony Guttmann~\cite{BrakGuttmann1990,FisherAuYang1979} developed methods
using algebraic and hypergeometric or holonomic functions to fit series data
starting from the early seventies. Of course, they named their techniques
differently, and, more importantly, they were primarily interested in
estimating \lq critical exponents\rq\ and \lq critical points\rq\ of the
function whose first few Taylor coefficients are given.

\section{Some Function Classes Suitable for Guessing}\label{sec:function-classes}
In this section we briefly present the function classes which are currently
explicitly covered by our package.  We want to stress however, that in many
cases it is easy to add other function classes, should the need arise.  (This
will become clear in Section~\ref{sec:background}.)

Throughout this article, $n\mapsto f(n)$ is the function we would like to
guess, and $f(x)=\sum_{n\ge0} f(n)x^n$ is its generating function.  The values
$f(n)$ are supposed to be elements of some field $\mathbb K$, usually the field
of rationals or rational functions. We alert the reader that the first value in
the given sequence always corresponds to the value $f(0)$.

\subsection{Guessing sequences $f(n)$}\label{sec:guess-fn}
\begin{description}
\item[\spad{guessRec}] finds recurrences of the form
  \begin{equation}\label{eq:Rec}
    p\big(f(n), f(n+1),\dots,f(n+k)\big)=0,
  \end{equation}
  where $p$ is a polynomial with coefficients in $\mathbb K[n]$.  For example,
\begin{verbatim}
guessRec [1,1,0,1,- 1,2,- 1,5,- 4,29,- 13,854,- 685]
\end{verbatim}
  yields $$[f(n): - f(n + 2) - f(n + 1) + f(n)^2 = 0,f(0)= 1,f(1)= 1].$$

  Note that, at least in the current implementation, we do not exclude
  solutions that do not determine the function $f$ completely.  For example,
  given a list containing only zeros and ones, one result will be
  \begin{equation*}
    [f(n): f(n)^2-f(n)=0,f(0)=\dots].
  \end{equation*}

\item[\spad{guessPRec}] only looks for recurrences with linear $p$, i.e., it
  recognises P-recursive sequences.  As an example,
\begin{verbatim}
guessPRec [0, 1, 0, -1/6, 0, 1/120, 0, -1/5040, 0, 1/362880, 
           0, -1/39916800, 0, 1/6227020800]
\end{verbatim}
  returns $$[f(n): (- n^2 - 3n - 2)f(n + 2) - f(n)= 0, f(0)= 0,f(1)= 1].$$

\item[\spad{guessRat}] finds rational functions.  For the sequence given in
  \eqref{eq:nsquared}, we find $n^2$ as likely solution.

\item[\spad{guessExpRat}] finds rational functions with an Abelian term, i.e.,
  \begin{equation*}
    f(n)=(a+bn)^n\frac{r(n)}{s(n)}
  \end{equation*}
  where $r$ and $s$ are polynomials.
\begin{verbatim}
guessExpRat [0,3,32,375,5184]
\end{verbatim}
  yields $$n(n+2)^n,$$ which could be interpreted, for example, as the number
  of labelled trees with one edge selected.

\item[\spad{guessBinRat}] finds rational functions with a binomial term, i.e.,
  \begin{equation*}
    f(n)=\binom{a+bn}{n}\frac{r(n)}{s(n)}
  \end{equation*}
  where $r$ and $s$ are polynomials.
\end{description}

Concerning $q$-analogues, \spad{guessRec(q)} finds recurrences of the
form~\eqref{eq:Rec}, where $p$ is a polynomial with coefficients in $\mathbb
K[q, q^n]$.  Similarly, we provide $q$-analogues for \spad{guessPRec} and
\spad{guessRat}.  For example, to guess a formula for
Sequence~\eqref{eq:binomial}, we enter\footnote{Because of a flaw in \FriCAS,
  one has to explicitly specify a list of options when using the $q$-versions
  of the guessing functions.  In the example above, we simply gave an empty
  list of options, and did thus not override any of the default options.}
\begin{verbatim}
guessRat(q)([1,1+q+q^2,(1+q+q^2)*(1+q^2),(1+q^2)*(1+q+q^2+q^3+q^4)], 
            [])
\end{verbatim}
and obtain as function
\begin{equation*}
  \frac{q^3q^{2n}+(-q^2-q)q^n+1}{q^3-q^2-q+1}.
\end{equation*}
Unfortunately the simplifying capabilities of \FriCAS\ are rather weak, so it
takes some extra work to simplify the above expression to
\begin{equation*}
  \frac{(1-q^{n+1})(1-q^{n+2})}{(1-q)(1-q^2)},
\end{equation*}
i.e., the $q$-binomial coefficient $\qbinom{n+2}{2}_q :=
\frac{[n+2]_q[n+1]_q}{[2]_q[1]_q}$, where
  $[n]_q := \frac{1-q^n}{1-q} = 1+q+\dots+q^{n-1}$.

Moreover, it is also possible to guess \lq mixed\rq\ recurrences, i.e., where
$p$ has coefficients in $\mathbb K[q, n, q^n]$, see the description of the
option \spad{maxMixedDegree} in Section~\ref{sec:options}.  For example,
\begin{verbatim}
guessPRec(q)([1,1,2*q^2,6*q^6,24*q^12,120*q^20,720*q^30,5040*q^42], 
             maxMixedDegree==2, homogeneous==true)
\end{verbatim}
returns $$[f(n): (n+1)f(n+1)q^{2n}-f(n+1)=0, f(0)=1].$$

The $q$-version of \spad{guessExpRat} recognises functions of the form
\begin{equation*}
  f(n)=(a+bq^n)^n\frac{r(q^n)}{s(q^n)},
\end{equation*}
$a$ and $b$ being in $\mathbb K[q]$ and $r$ and $s$ polynomials with
coefficients in $\mathbb K[q]$.  For Sequence~\eqref{eq:qExpRat}, we enter
\begin{verbatim}
guessExpRat(q)([(1-2*q)/(1-q),1-2*q,(1-q)*(1-2*q)^3,
                (1-q)^2*(1-2*q)*(1-2*q-2*q^2)^3], [])
\end{verbatim}
to obtain
\begin{equation*}
  \frac{2q-1}{q-1}(2q^n-3q+1)^n.
\end{equation*}
Another example would be Nicholas Loehr's $q$-analogue $[n+1]_q^{n-1}$ of
Cayley's formula.

Finally, \spad{guessBinRat(q)} tries to fit the
given terms to 
\begin{equation*}
  f(n)=\qbinom{a+bn}{n}_q\frac{r(q^n)}{s(q^n)},
\end{equation*}
where $\qbinom{n}{m}_q=\prod_{i=1}^m\frac{1-q^{n-i+1}}{1-q^i}$.

\subsection{Guessing Series $f(x)$}\label{sec:guess-fx}

\begin{description}
\item[\spad{guessADE}] finds an algebraic differential equation for $f(x)$,
  i.e., an equation of the form
  \begin{equation}
    \label{eq:ADE}
    p\left(f(x), f^\prime(x),\dots,f^{(k)}(x)\right)=0,
  \end{equation}
  where $p$ is a polynomial with coefficients in $\mathbb K[x]$. A typical
  example is $\sum n^n\frac{x^n}{n!}$:
\begin{verbatim}
guessADE [1,1,2,9/2,32/3,625/24,324/5,117649/720,131072/315,
          4782969/4480]
\end{verbatim}
  returns
\begin{equation*}
  [[x^n ]f(x): - xf'(x) + f(x)^3  - f(x)^2 = 0,f(0)= 1,f'(0)= 1].
\end{equation*}

Maybe more interesting, we obtain also a differential equation for the
exponential generating function with coefficients of the form covered by
\spad{guessExpRat}:
\begin{verbatim}
guessADE([(a*n+b)^n/factorial(n) for n in 0..32], 
         maxPower==3, maxDerivative==3, homogeneous==true)
\end{verbatim}
However, this equation is already quite big:
\begin{multline*}
  4 b^2 (a+b)^2 f(x)^2 f''(x) + 3 a b^2 (a+b) x f(x)^2 f'''(x)
  - 4b^2(2a+b)^2 f(x) f'(x)^2\\
  - a(a^3+a^2b+19ab^2+3b^3)xf(x)f'(x)f''(x)
  - 3a^3(a-3b) f(x) f'(x) f'''(x)\\
  + 5a^3(a-3b)x^2 f(x)f''(x)^2 + 4a^4 x f'(x)^3=0.
\end{multline*}
We stress that we did not try to prove this equation -- it remains a guess,
even though we checked the first few hundred terms.

Another interesting example is given by the generating function for the
chromatic polynomials of rooted triangulations, as found by William
Tutte~\cite{MR1336864}.  Or, as a test case, to guess a differential equation
for Jacobi's $\theta$-function $1+2\sum z^{n^2}$ a list of the first $3600$
terms,
\begin{verbatim}
guessADE(l, maxPower==14, maxDerivative==3, maxDegree==6)
\end{verbatim}
and a little patience (roughly ten minutes on an AMD Opteron processor)
suffice.  In fact, according to Don Zagier~\cite[Section~5.1,
Proposition~15]{MR2409678} already Ramanujan knew that every modular and every
quasi-modular form on $\Gamma_1$ satisfies a third order algebraic differential
equation.

\item[\spad{guessHolo}] only looks for equations of the form~\eqref{eq:ADE}
  with linear $p$, that is, it recognises holonomic or differentially-finite
  functions.  It is well known that the class of holonomic functions coincides
  with the class of functions having P-recursive Taylor coefficients.  However,
  the number of terms necessary to find the differential equation often differs
  greatly from the number of terms necessary to find the recurrence.  Returning
  to the example given for \spad{guessPRec}, we find that already
\begin{verbatim}
guessHolo [0,1,0,-1/6,0,1/120]
\end{verbatim}
  returns
  \begin{equation*}
    [[x^n ]f(x):  - f''(x) - f(x)= 0,f(0)= 0,f' (0)= 1].
  \end{equation*}
  Moreover, now we immediately recognise the coefficients as being those of the
  sine function.

\item[\spad{guessAlg}] looks for an algebraic equation satisfied by $f(x)$,
  i.e., an equation of the form
  \begin{equation*}
    p\left(f(x)\right)=0,
  \end{equation*}
  the prime example being given by the Catalan numbers
\begin{verbatim}
guessAlg [1,1,2,5,14,42]
\end{verbatim}
  which yields
\begin{equation*}
  [[x^n ]f(x): x f(x)^2  - f(x) + 1= 0,f(0)= 1].
\end{equation*}

\item[\spad{guessPade}] recognises rational generating functions or,
  equivalently, recurrences with constant coefficients.  For the Fibonacci
  sequence given in \eqref{eq:Fibonacci}, we find as likely solution
\begin{equation*}
  [[x^n ]f(x): (x^2  + x - 1)f(x) + 1= 0].
\end{equation*}

\item[\spad{guessFE}] finds \lq Mahler-type\rq\ functional equations for $f(x)$
  (see for example \cite{Mahler1976}), i.e., equations of the form
  \begin{equation}
    \label{eq:FE}
    p\left(f(x), f(x^2),\dots,f(x^k)\right)=0,
  \end{equation}
  where $p$ is a polynomial with coefficients in $\mathbb K[x]$.  A typical
  example is the number of unlabelled rooted binary trees:
\begin{verbatim}
guessFE [0,1,1,1,2,3,6,11,23]
\end{verbatim}
  which returns
\begin{equation*}
  [[x^n]f(x): f(x^2) + f(x)^2  - 2f(x) + 2x= 0, f(x)= x + x^2+ x^3+2x^4+O(x^5)].
\end{equation*}

Browsing \emph{the online encyclopedia of integer sequences}, we discovered
another rather surprising functional equation: consider the sequence
\OEIS{A118006} of binary words $w_n$ defined by $w_1="01"$ and $w_{n+1}=concat
[w_n, w_n, reverse(w_n)]$. Then
\begin{verbatim}
guessFE w 4
\end{verbatim}
indicates that the limiting word $w_\infty$ satisfies
\begin{equation*}
  (x-1)(x^2-x+1)(x^2+x+1)((x^2+x+1)\left(f(x)-x^2 f(x^3)\right) + x(x^2+1)^2=0.
\end{equation*}
Again, we did not try to prove this equation but only checked the first few
hundred terms.
\end{description}

For \spad{guessADE} and \spad{guessHolo} we provide $q$-analogues, replacing
differentiation with $q$-dilation: \spad{guessADE(q)} finds differential
equations of the form
\begin{equation}
  \label{eq:qADE}
  p\left(f(x), f(qx),\dots,f(q^kx)\right)=0,
\end{equation}
where $p$ is a polynomial with coefficients in $\mathbb K[q, x]$.  Generating
functions satisfying such $q$-equations frequently occur in the enumeration of
polyominoes and the study of orthogonal polynomials.  As an example, we can
recover the $q$-algebraic differential equation for the generating function of
bar polyominoes by horizontal perimeter -- marked by $x$, vertical perimeter --
marked by $y$ and area -- marked by $q$, as given by Richard Brak and Thomas
Prellberg~\cite{BrakPrellberg1995}.  We enter
\begin{verbatim}
guessADE(q)(l, maxDerivative==1, maxPower==2, maxDegree==1)
\end{verbatim}
where \spad{l} are the first eleven coefficients of the series in $x$:
\begin{verbatim}
l := [0, q*y/(1-q*y), q^2*y*(1+q*y)/(1-q*y)/(1-q^2*y), ...]
\end{verbatim}
The solver immediately finds the solution
\begin{multline*}
[x^n] f(x): (qxf(x)+(qx+1)y)f(qx)+(qx-1)f(x)+qxy=0,\\
f(0)=0, f'(0)=\frac{qy}{1-qy},\dots],
\end{multline*}
it then takes a few seconds to verify it.

\subsection{Operators}\label{sec:operators}

The observation made by Christian Krattenthaler before writing his program
\Rate~\cite{RATE} is the following: it occurs frequently that although a
sequence of numbers is not generated by a rational function, the sequence of
successive quotients is.

We slightly extend upon this idea, and apply recursively one or both of the two
following operators:
\begin{description}
\item[\spad{guessSum} - $\Delta_n$] the differencing operator, transforming
  $f(n)$ into $f(n)-f(n-1)$.
\item[\spad{guessProduct} - $Q_n$] the operator that transforms $f(n)$ into
  $f(n)/f(n-1)$.
\end{description}

For example, to guess a formula for Sequence~~\eqref{eq:sumfact}, we enter
\begin{verbatim}
guess([0, 1, 3, 9, 33], [guessRat], [guessSum, guessProduct]).
\end{verbatim}
The second argument to \spad{guess} indicates which of the functions of the
previous section to apply to each of the generated sequence, while the third
argument indicates which operators to use to generate new sequences.

The package will then respond with
\begin{equation*}
  \sum_{s=0}^{n-1}\prod_{p=0}^{s-1} (p+2),
\end{equation*}
i.e., the sum of the first factorials.

In the case where only the operator $Q_n$ is applied, our package is directly
comparable to \Rate. In this case the standard example is the number of
alternating sign matrices
\begin{verbatim}
guess [1, 1, 2, 7, 42, 429, 7436, 218348]
\end{verbatim}
which yields
\begin{equation*}
  \prod_{k=0}^{n-1}\prod_{l=0}^{k-1}\frac{27l^2+54l+24}{16l^2+32l+12}
  =\prod_{k=0}^{n-1}\prod_{l=0}^{k-1}\frac{3(3l+2)(3l+4)}{4(2l+1)(2l+3)}.
\end{equation*}

\subsection{Closure properties and zero test}
\label{sec:remarks}

Part of what makes a class of functions interesting are its closure properties,
summarised in the table below for some classes of functions.  Apart from the
theoretical point of view, it is also good to know that the computer can guess
an equation for $f(n)$ if it can do so for $f(n)+1$.

However, one has to keep in mind that even simple transformations may increase
the number of terms necessary to successfully guess an equation dramatically.
For example, consider the (exponential) generating function for the Bell
numbers $B_n$, counting the number of partitions of $\{1,2,\dots,n\}$, which is
\begin{equation*}
  B(x)=\sum_{n\ge0} B_n \frac{x^n}{n!} = e^{e^x-1}.
\end{equation*}
This series is not holonomic, but it satisfies the simple algebraic
differential equation $B'' B - (B')^2 - B'B =0$, and the first thirteen terms
suffice to find it.  By contrast, it takes $36$ terms to guess the shifted
series $(e^{e^x-1}-1)/x$.

In the same spirit, note that without specifying the search space any further,
already six terms are enough to guess a functional equation for the number of
unlabelled rooted binary trees.  On the other hand, we need at least $42$ terms
to guess an equation for the square of their generating function.

This phenomenon also explains why Christian Krattenthaler's program \Rate\ is
so useful: of course there is also an algebraic recurrence for the number of
alternating sign matrices, namely
\begin{equation*}
  (- 16n^2  - 32n - 12)f(n)f(n + 2) + (27n^2  + 54n + 24)f(n + 1)^2 = 0,
\end{equation*}
but we need $35$ terms to guess it instead of eight.  (Instead of looking for a
formula having $k$ nested products, we could also use the options
\spad{Somos==true, maxShift==k, homogeneous==2\^{}(k-1)}, but this only works
well for $k$ less than 4.)
%

\begin{table}[h]
  \centering
\def\y{\checkmark}
\def\a{alg.}
\def\ca{$\cdot$}
\def\iv{$(.)^{-1}$}
\def\co{$\circ$}
\def\ci{$(.)^{(-1)}$}
\def\In{$\int$}
\def\ha{$\odot$}
\def\sh{S}
\begin{tabular}{|l||c|c|c|c|c|c|c|c|c|}\hline
 type of equation     &$+$&\ca&\iv&\co&\ci&$D$&\In&\ha&\sh\\\hline\hline
                Pad\'e&\y &\y &\y &\y & - &\y & - &\y &\y\\\hline
             algebraic&\y &\y &\y &\y &\y &\y & - & - &\y\\\hline
   linear differential&\y &\y & - &\a & - &\y &\y &\y &\y\\\hline
algebraic differential&\y &\y &\y &\y &\y &\y &\y & - &\y\\\hline
\end{tabular}
\medskip
\caption{closure properties. (\iv:~multiplicative inverse, \co:~composition,
  \ci:~compositional inverse, \In:~definite integration, \ha:~Hadamard product,
  \sh:~shift, \a:~algebraic substitution; inverse and substitution only apply when the result is again
  a formal power series.} 
\end{table}

Proofs for the closure properties of
rational, algebraic and linear differential equations can be found, for
example, in Richard Stanley's book Enumerative Combinatorics~2~\cite{EC2} or
his article on differentiably finite series~\cite{Stanley1980}.  
For algebraic differential equations, proofs
were given by Alexander Ostrowski in~\cite{Ostrowski1920}, see
also~\cite{Klazar2003}.  Slightly weaker closure
properties hold for the $q$ case.  In particular, $q$-holonomic series are only
closed under the substitution $x\mapsto x^k, k\in\mathbb N$, see for example
\cite{RISC3145}.

Algebraic recurrence relations seem to satisfy no
interesting closure properties: for example, take any sequence that
does not satisfy an algebraic recurrence relation, and write it as
the sum of two sequences, one with odd terms zero, the other with
even terms zero.  Both summands are solutions of $f(n)f(n+1) = 0$.
However, a related class, so called \lq
admissible recurrences\rq\ was studied by Manuel Kauers~\cite{Kauers2005} and
has been shown to enjoy many closure properties.  

Similarly, we are not aware of any results concerning closure
properties of Mahler-type functional equations as defined in our
paper.  However, linear equations $p\left(f(x),
f(x^r),f(x^{r^2})\dots,f(x^{r^k})\right)=0$ for fixed $r$ were shown
by Phillippe Dumas~\cite{MR1346304} to give rise to nice closure
properties.  This extends to the non-linear case.

One very important property of these classes is the availability of a
\Dfn{zero-test}, i.e., an algorithm that will decide whether any given equation
(together with sufficiently many initial values) has only the zero solution.
For linear differential equations this is folklore, for algebraic differential
equations an algorithm was proposed for example by Joris van der
Hoeven~\cite{vanDerHoeven}.  In many cases, such a test allows to verify
conjectured identities automatically, as exercised for example by \GFUN.

\section{Options}\label{sec:options}
To give you the maximum flexibility in guessing a formula for your favourite
sequence, we provide options that modify the behaviour of the functions as
described in Section~\ref{sec:function-classes}. The options are appended,
separated by commas, to the guessing function in the form \spad{option==value}.
See below for some examples.

\begin{description}
\item[\spad{maxDerivative}, \spad{maxShift}] specify the maximum derivative in
  an algebraic differential equation, or, in a recurrence relation, the maximum
  shift.  Setting the option to \spad{arbitrary} specifies that the maximum
  derivative -- the maximum shift -- may be arbitrary, which is the default.

\item[\spad{maxPower}] specifies the maximum total degree in an algebraic
  differential equation or recurrence: for example, the degree of $(f'')^3 f'$
  is $4$. Setting the option to \spad{arbitrary} specifies that the maximum
  total degree may be arbitrary, which is the default.  

\item[\spad{homogeneous}] specifies whether the search space should be
  restricted to homogeneous algebraic differential equations or homogeneous
  recurrences, i.e., the case where the polynomial $p$ in
  Equation~\eqref{eq:Rec} and Equation~\eqref{eq:ADE} is homogeneous.  By
  default, it is set to \spad{false}.  Setting it to a positive integer, only
  homogeneous polynomials $p$ of this degree are tried.  Setting it to
  \spad{true}, all homogeneous polynomials $p$ up to total degree
  \spad{maxPower} are tried.

\item[\spad{Somos}] specifies whether the search space should be restricted to
  algebraic differential equations where the sum of differentials is constant.
  Similarly, when guessing recurrences, \spad{Somos} insists that the sum of
  shifts is constant.  By default it is set to \spad{false}.  Setting it to a
  positive integer, the sum of differentials or shifts must be equal to this
  number.  Setting it to \spad{true} is equivalent to invoking the guesser with
  \spad{Somos==2}, \spad{Somos==3}, \dots, \spad{Somos==d}, where \spad{d} is
  the specified \spad{maxDerivative} (or \spad{maxShift}) times \spad{maxPower}
  or \spad{homogeneous}.

\item[\spad{maxDegree}] specifies the maximum degree of the coefficient
  polynomials in an algebraic differential equation, a Mahler-type functional
  equation or a recurrence with polynomial coefficients.  For rational
  functions with an exponential term, \spad{maxDegree} bounds the degree of the
  denominator polynomial.  The default value of \spad{maxDegree} is
  \spad{arbitrary}.

\item[\spad{allDegrees}] specifies whether all possibilities of the degree
  vector -- taking into account \spad{maxDegree} -- should be tried.  The
  default is \spad{true} for \spad{guessPade} and \spad{guessRat} and
  \spad{false} for all other functions.

\item[\spad{maxMixedDegree}] allows guessing of mixed $q$-recurrences. Its
  value determines the maximum degree of $q^n$ in the coefficients, default
  being zero.

\item[\spad{maxLevel}] specifies how many levels of recursion are tried when
  applying operators as described in Section~\ref{sec:operators}.  Note that,
  applying either of the two operators results in a sequence which is by one
  shorter than the original sequence.  Therefore, in case both \spad{guessSum}
  and \spad{guessProduct} are specified, the number of times a guessing
  algorithm from the given list of functions is applied is roughly $2^n$, where
  $n$ is the number of terms in the given sequence.  Thus, especially when the
  list of terms is long, it is important to set \spad{maxLevel} to a low value.

  Still, the default value is \spad{arbitrary}, which means that the number of
  levels is only restricted by the number of terms given in the sequence.

\item[\spad{safety}] specifies, as explained in detail in
  Section~\ref{sec:background} and Section~\ref{sec:safety} the number of
  additional equations a solution has to satisfy.  The default setting is $1$.

  Experiments indicate that, the larger the class of functions covered, the
  larger one should set \spad{safety}.  Moreover, when the sequence contains
  many zeros, higher settings of safety are appropriate.  For all algorithms we
  recommend to set \spad{safety} higher than the number of trailing zeros.  The
  reason is best illustrated by an example:
\begin{verbatim}
guessPade([a,b,c,0])
\end{verbatim}
  returns
  \begin{equation*}
    [[x^n]c x^2  + b x + a].
  \end{equation*}
  In other words, if the sequence has a trailing zero, \spad{guessPade}
  trivially finds a solution.  A few experiments and a moment's thought will
  reveal that the other algorithms behave similarly.

\item[\spad{check}] determines whether we want to check the solutions
  returned by the modular solver using a deterministic check, or
  whether we content ourselves with a (rather weak) Monte-Carlo type
  check, or skip checking entirely, the default value being
  \spad{deterministic}.

\item[\spad{checkExtraValues}] specifies whether we want to return only those solutions
  that fit the given data perfectly.  With \spad{checkExtraValues==false}, the complete
  basis of the solution space is returned, see Section~\ref{sec:safety}.  The default value
  is \spad{true}.

\item[\spad{one}] specifies whether the guessing function should return as soon
  as at least one solution is found. By default, this option is set to
  \spad{true}.

\item[\spad{indexName}, \spad{variableName}, \spad{functionName}] specify
  symbols to be used for the output. The defaults are \spad{n}, \spad{x} and
  \spad{f} respectively.

\item[\spad{debug}] specifies whether information about progress should be
  reported.
\end{description}

\section{Rational Interpolation}\label{sec:background}

The underlying idea of all guessing software is to fit the given data to a
\Dfn{model}.  For example, a formula for Sequence~\eqref{eq:nsquared}, is
almost trivial to guess: it seems obvious that it is $n^2$.  A natural model to
check is that the sequence in question is generated by a polynomial -- we
simply apply polynomial interpolation.  Given a list of four terms -- $0,1,4,9$
in our example -- we should expect that we need a polynomial of degree three to
interpolate.  Since the actual degree is lower, that is, the interpolating
polynomial is overdetermined by the data, it is reasonable to accept $n^2$ as a
good guess.

Generalising to \Dfn{Hermite-Pad\'e interpolation}, we can cover most models
described in Section~\ref{sec:guess-fn}:
\begin{prbSequence}
  Let $\Mat f = [f^{(1)}(n),\dots,f^{(m)}(n)]$ be a vector of (truncated)
  sequences over some integral domain, and $\Mat n = [n^{(1)},\dots,n^{(m)}]$ a vector of
  non-negative integers, serving as degree bounds.  Let
  $\sigma\geq0$.  Determine a polynomial vector
  $\Mat p = [p^{(1)}(n),\dots,p^{(m)}(n)]$ with $\deg p^{(l)}(n) < n^{(l)}$
  such that
  \begin{equation}
    \label{eq:orderSequence}
    p^{(1)}(n)\cdot f^{(1)}(n) + \dots + p^{(m)}(n)\cdot f^{(m)}(n)=0\text{ for }
    0\leq n<\sigma.
  \end{equation}
\end{prbSequence}

Note that, by equating coefficients, this problems can be reduced to solving an
appropriate \emph{linear} system of equations with  $n^{(1)}+\dots+n^{(m)}-1$
unknowns, namely the coefficients of the polynomials
$p^{(1)}(n),\dots,p^{(m)}(n)$, up to normalisation.  Thus, we will in fact
determine a basis of the space of solutions.
However, instead of using, for example, naive Gaussian elimination, we will
take advantage of the special structure of these linear systems to achieve
better performance.  To illustrate, we would like to be able to solve systems
where $n^{(1)}+\dots+n^{(m)}$ is as large as ten thousand.

Setting $\sigma=n^{(1)}+n^{(2)}-1$ and 
$\Mat f=[(1,1,\dots,1),(f_0,f_1,\dots,f_{\sigma-1})]$ we would recover ordinary
rational interpolation.  However, to have more confidence in the \lq
guessed\rq\ formula, we use
$\sigma=n^{(1)}+\dots+n^{(m)}-1+\spad{safety}$ instead.

More generally, to guess algebraic recurrences we
consider the (infinite) sequence of monomials in the \lq variables\rq\ $f(n),
f(n+1), f(n+2), \dots$
$$\big(\prod_i f(n+\lambda_i-1)\big)_\lambda,$$
where $\lambda=(\lambda_1,\lambda_2,\dots)$ runs over the integer partitions in
lexicographic order:
\begin{equation*}
  1, f(n), f(n)^2, f(n+1), f(n)^3, f(n) f(n+1), f(n+2), f(n)^4, f(n)^2 f(n+1),\dots
\end{equation*}
Then, for each $m\geq 2$ we solve the rational interpolation problem with $\Mat
f$ given by the first $m$ entries of this sequence,
and $\Mat n$ such that the number of unknowns
$n^{(1)}+\dots+n^{(m)}-1$ in the corresponding linear system plus the specified
value of \spad{safety} equals the number of equations $\sigma$.

In the formulation of the rational interpolation problem above, the sequence of
evaluation points was chosen as $0,1,2,\dots$, but it is straightforward to
generalise to arbitrary evaluation points.  Doing so, we can also find
$q$-recurrences, by pretending that $f$ is given at the points
$q^0,q^1,q^2,\dots$ instead.

To deal with the models described in Section~\ref{sec:guess-fx}, we need to
solve another variant of the rational interpolation problem:
\begin{prbSeries}
  Let $\Mat f = [f^{(1)}(x),\dots,f^{(m)}(x)]$ be a vector of (truncated) power
  series over some integral domain, and $\Mat n = [n^{(1)},\dots,n^{(m)}]$ a vector of
  non-negative integers, serving as degree bounds.  Let
  $\sigma\geq0$. Determine a polynomial vector
  $\Mat p = [p^{(1)}(x),\dots,p^{(m)}(x)]$ with $\deg p^{(l)}(x) < n^{(l)}$
  such that
  \begin{equation}
    \label{eq:orderSeries}
    \ord\left(\Mat p\cdot\Mat f\right) = \ord\left(p^{(1)}(x)\cdot f^{(1)}(x) +
      \dots + p^{(m)}(x)\cdot f^{(m)}(x)\right) \geq\sigma.
  \end{equation}
\end{prbSeries}

In this case, setting $\sigma=n^{(1)}+n^{(2)}-1$ and  
$\Mat f=[1,f(x)]$, where $f(x)$ is the truncated power
series with the given values as Taylor coefficients, we recover Pad\'e
approximation.  This allows us to \lq guess\rq\ sequences that are Taylor
coefficients of rational generating functions. 

To guess algebraic differential equations, we consider the sequence of
monomials $\left(\prod_i \D^{\lambda_i-1}f(x)\right)_\lambda$, where $\D$ is
the differentiation operator and $\lambda=(\lambda_1,\lambda_2,\dots)$ runs
over the integer partitions in lexicographic order as before:
\begin{equation*}
  1, f(x), f(x)^2, f'(x), f(x)^3, f(x)f'(x), f''(x), f(x)^4, f(x)^2f'(x),\dots
\end{equation*}

To guess $q$-algebraic differential equations, we just replace the usual
differentiation operator with $q$-dilation: $\D f(x) := f(qx)$.  Finally,
\spad{guessFE} uses the sequence of monomials $\left(\prod_i
  f(x^{\lambda_i})\right)_\lambda$.

For the present package, we originally implemented a fraction free algorithm
proposed in 2000 by Bernhard Beckermann and George
Labahn~\cite{BeckermannLabahn2000}, which at the time proved much faster than
what \GFUN\ had.  However, during the refereeing process it became clear that a
modular approach would be even more efficient.  This was first pointed out by
Manuel Kauers, and independently by Alin Bostan and Bruno Salvy in private
communication.  Consequently, we decided to follow this approach and
implemented a modular version of an \emph{older} algorithm from 1994, also by
Bernhard Beckermann and George Labahn~\cite{BeckermannLabahn1994}, when the
coefficients are rational numbers or rational functions with integer
coefficients.  This turned out to be very fruitful, although quite
labour-some. For other coefficient domains we still use the fraction
free algorithm, although we plan to extend the modular approach to
allow algebraic numbers as coefficients as soon as possible.

We would like to stress that meanwhile most of the packages mentioned in the
introduction use modular techniques, however using other algorithms for solving
over a prime field.  According to Bruno Salvy, \GFUN\ now uses an algorithm
introduced in 1997, again by Bernhard Beckermann and George
Labahn~\cite{BeckermannLabahn1997}.  Manuel Kauers package \GuessKauers\ uses
the solver provided by \MMA, it is thus unclear which algorithm is used.

Still, our package outperforms the other freely available packages,
for many configurations of degree bounds and size of the vector $\Mat
f$, (see Section~\ref{sec:performance}), as well as -- for univariate
sequences -- the range of formulae that can be guessed.

We also implemented specialised algorithms to test whether the
$n$\textsuperscript{th} term of the sequence is given by a formula of the form
\begin{equation}\label{eq:abel}
  n\mapsto (a+bn)^n\frac{r(n)}{s(n)}\quad\text{or}\quad
  n\mapsto\binom{a+bn}{n}\frac{r(n)}{s(n)},
\end{equation}
for some $a$ and $b$ and polynomials $r$ and $s$.  Unfortunately, we could not
avoid solving \emph{non-linear} equations in this case.  Even after exploiting
some surprising coincidences that reduce the size of the arising equations the
performance of this algorithm is disappointing: already eight or nine terms,
i.e., degree two in $r$ and $s$ pose a challenge, even over a finite
field.\footnote{Meanwhile, it seems that we have found a suitable approach, but
  due to time constraints we cannot describe it in this article.}

\section{Safety}\label{sec:safety}
How can we \lq know\rq\ that a formula discovered via interpolation is
appropriate?  At first glance, the answer is quite simple: we use all but the
last few terms of the sequence to derive the formula.  After this, the last
terms are compared with the values predicted by the polynomial.  If they
coincide, we can be confident that the guessed formula is correct.

In the case of the rational interpolation problem we get the same set
of accepted solutions when we use all values, but keep lower
degree bounds.  We use this approach as it is more efficient
than actually computing \lq bad\rq\ solutions and rejecting them
later, although there is a subtle interaction with an extra check that
we perform.

Very recently, Alin Bostan and Manuel 
Kauers~\cite[Section~2.4]{BostanKauers2009} described in some detail various
other possibilities of checking whether a guessed formula is likely to be \lq
correct\rq, the method we just outlined being clearly the most practical.
Unfortunately, it turns out that this method is problematic in certain
situations.  In this section we explain why.

First of all, we cannot expect that all elements of the solution
space of the rational interpolation problem \lq interpolate\rq\ the
given data in the following sense: consider the truncated power
series $f(x)=1+x^6+O(x^7)$, and let $\Mat f = [1, f(x), f'(x)]$.
Setting the vector of degree bounds $\Mat n = [2,2,2]$ and $\sigma=6$
(note that we \lq loose\rq\ one term because of differentiation,
so we have $6$ equations in our linear system), 
rational interpolation yields the basis $[(1, -1, 0), (0, 0, x)]$.
Thus, the general solution to the rational interpolation problem with
the given constraints is
$$(\alpha+\beta x)\left(1-f\left(x\right)\right) + \gamma x f'(x)=0,$$
$\alpha$, $\beta$ and $\gamma$ being elements from the coefficient
field.

Apparently, none of the two basis vectors actually interpolates all
given values: $1-f\left(x\right)=-x^6+O(x^7)$, and $x f'(x) =
6x^6+O(x^7)$.  One might be tempted to simply discard
non-interpolating basis vectors (which we do when
\spad{checkExtraValues} is true), but doing so we risk loosing \lq
good\rq\ solutions, too: 
$$(6\gamma + \beta x)\left(1-f\left(x\right)\right) + \gamma x f'(x) = O(x^7)$$
interpolates just fine
for any $\beta$ and $\gamma$.  In particular, the set of interpolating solutions
is not a vector space.

An uncomfortable consequence of the above is as follows: we provide an option
\spad{maxDegree} that allows the user to specify the maximum degree of the
coefficient polynomials, see Section~\ref{sec:options}.  When set to some
integer value $d$, we (essentially) do not compute solution spaces of
configurations $\Mat f$ with $(d+1)\size{\Mat f}$ being less than the number of
values provided.  Suppose now that we find an interpolating solution without
setting \spad{maxDegree}, and that the maximal coefficient degree of this
solution happens to be $d$.  Then it may be the case that setting
\spad{maxDegree==d} instead yields no result, because all basis elements are
discarded.  Similarly, one might expect that increasing
both \spad{safety} and the number of values by one does not yield
more solutions.  But at lower \spad{safety}
our check may reject all basis elements, while at
higher \spad{safety} the basis may contain an interpolating solution.

%

A possible way to resolve this dilemma might be to reject solution spaces that
are not one-dimensional.  However, when pursuing this idea, another difficulty
surfaces: namely, it is not completely trivial to decide whether two solutions
are \emph{really} different.  For example, consider $\Mat f = [1, f(x),
f'(x)]$, and suppose that $f(x)$ is in fact a polynomial $p(x)$.  Then the
interpolation routine will not only find the solution $[p(x), 1, 0]$, but also
$[p'(x), 0, 1]$. 
More generally, it is well known that one often needs more
coefficients to determine the minimal order equation than to find a
solution of higher order.  Thus, if we have enough values to guess
the minimal order equation then the problem is easy.  But otherwise
we will either find multiples of the minimal equation, or some
parasitic solutions.

This problem can
be remedied, at least for linear and also algebraic differential equations: in
the linear case, we could simply compute a greatest common right divisor of the
given equations, whereas in the algebraic case we could apply Ritt
elimination.

Still, there is again some danger that \lq good\rq\ solutions are lost: for
example, if a sequence is non-zero only at very few indices $n_1,
n_2,\dots,n_k$, then the interpolation algorithm will not only find the \lq
good\rq\ solution, but also $(n-n_1)(n-n_2)\dots(n-n_k)f(n)=0$, and the
greatest common right divisor of the two will be trivial.

We admit that so far we were unable to find a completely satisfying solution to
this problem.  In the meantime, we provide options (in particular
\spad{checkExtraValues},
and \spad{one}, see Section~\ref{sec:options}) that let the user decide.

\section{Rational Reconstruction}
\label{sec:rational-reconstruction}

As already mentioned in Section~\ref{sec:background}, our solver uses a modular
technique: instead of solving the rational interpolation problem over the
integers, we solve the problem over several machine size prime fields and use
Chinese remaindering to obtain integer solutions.  In the same spirit, given
coefficients that are rational functions, we evaluate them at several random
points, solve the simpler problems and use rational reconstruction to obtain
polynomial solutions.

There are two different ways in which this plan can fail: it may happen that
the solution of the problem in the prime field is not an image of the solution
of the problem in the original ring.  In
Corollary~\ref{cor:bad-reduction} we will see that there are ways to discover
such \lq bad reductions\rq, provided we have at least one \lq good
reduction\rq.  However, we cannot a priori exclude the possibility that all
reductions are bad.

Moreover, it may occur that due to an unfortunate choice of evaluation points
we obtain wrong solutions -- usually, when we have too few evaluation points we
get no solution, but it may happen that we construct one that is actually
wrong.

Therefore, the solution returned from the core solver is only probably correct,
and we need to check it before returning it.  Thus, the main loop of the solver
in pseudocode is:
\begin{verbatim}
repeat
   sol := do_solve(data, inner_call? == false)
   if check(sol) then return sol
\end{verbatim}
where \verb/do_solve/ produces a probably correct solution and \verb/check/
verifies correctness.  So far we did not encounter a case where the check failed
-- the solver is designed in such a way that probability of wrong answer is
very low. Therefore, instead of looping, we print an error message and fail.

\verb/do_solve/ is a Brown-style~\cite{MR0307450} routine similar to
Subroutine~M and Subroutine~P in the gcd algorithm of Mark van Hoeij and
Michael Monagan~\cite{MR2126957}, which in pseudocode looks as follows:
\begin{verbatim}
do_solve(data) ==
    if R = Z_p then return solve_over_Z_p(data)
    bad_count := 0
    good_count := 0
    sol := empty()
    repeat
        modulus := choose_modulus()
        if inner_call? then
            new_data := eval(data, modulus)
        else
            new_data := reduce(data, modulus)
        new_sol := do_solve(new_data, inner_call? == true)

        if new_sol = "no_solution" then return "no_solution"
        reduction_status := check_reduction(new_sol)
        if new_sol = "failed" or reduction_status = "bad" then
            bad_count := bad_count + 1
            if inner_call? and bad_count > good_count + 2 then 
                return "failed"
        else
            good_count := good_count + 1
            if reduction_status = "all_bad" then sol := empty()
            sol := chinese_remainder(sol, new_sol, modulus)
            rr := rational_reconstruction(sol)
            if not rr = "failed" then return rr
\end{verbatim}
In contrast to Mark van Hoeij and Michael Monagan we present this algorithm as
a single routine, to stress that the processing is generic: in the outer level,
when \verb/inner_call?/ is false, 
\verb/choose_modulus/ chooses machine size primes, in the inner level it
chooses random evaluation points.

The routine \verb/check_reduction/ applies Corollary~\ref{cor:bad-reduction} to
\verb/new_sol/ -- we keep the necessary information about previously
obtained solutions in a global variable, namely the minimal dimension
of the solution space and the minimal values of the critical indices.
(Of course, when there are no previous solutions, \verb/new_sol/ is
automatically treated as \lq good\rq.)  Also, we discard solutions
with leading exponent being smaller than the leading exponent of previous solutions --
this is necessary to ensure correct normalisation
after rational reconstruction and also avoids the problem of \lq bad
content\rq, see~\cite{MR2126957}.

The counters \verb/bad_count/ and
\verb/good_count/ are used to detect cases where we encounter bad reduction
already at the outer level -- we copied the method used in \cite{MR2126957}.
Without this test the solver would spend a lot of time reconstructing solutions
which would then be discarded at outer level.

The variable \verb/sol/ contains homomorphic images of the bases of the
solution spaces constructed in the current stage, and
\verb/rational_reconstruction(sol)/ tries to find the basis in the original
ring.  In the following, we indicate which tricks we decided to implement to
make the procedure efficient enough.
Namely, for polynomials we use naive quadratic multiplication, gcd and Chinese
remaindering (Lagrange interpolation) routines.  Also our rational
reconstruction implementation uses a simple quadratic method.
However, we save time by trying 
rational reconstruction not in every step but only after an interval: 
for polynomials we use a quadratically growing sequence of points, while for
integers we switch to a big step (currently $100$) after passing a threshold
(currently $200$).

More precisely, to reconstruct a solution in characteristic zero
given modular images $m_1$, $m_2$, \dots, $m_n$ and moduli $p_1$,
$p_2$, \dots, $p_n$ we need to compute $M$ such that $M = m_i\mod
p_i$, and then apply rational reconstruction to $M$.  However,
instead of computing $M$ directly we first compute intermediate
solutions $M_j$ such that $M_j = m_i\mod p_i$ for $i\in\{100 j,\dots,
100 j + 99\}$, updating incrementally.  Whenever we have finished
such a block of 100 primes and computed $M_j$, we update $M$ such
that $M = M_j\mod P_j$ where $P_j = p_{100 j} p_{100 j+1}\dots p_{100
  j+99}$, and then apply rational reconstruction to $M$.  This scheme
makes more efficient use of bignum routines: we perform most
operations on relatively small bignums, and fewer operations on big
bignums.

There is one more improvement to rational reconstruction which 
first appeared in NTL~\cite{Shoup2005} (see also Section~3.1 of the
description of the IML~\cite{MR2280534}, which is where we learned
from the trick): when reconstructing the vector of 
rational coefficients, we incrementally compute the common denominator
of the coefficients already reconstructed, and impose it on subsequent terms.
Since we are looking for an integer (or polynomial) solution and fractions
appear only due to normalisation, it is natural to expect
all terms to have the same or very similar denominators.
Often, multiplying the next term by the common denominator computed
so far, it turns out that the product is already acceptable as
rational reconstruction of this coefficient.

\section{Implementation aspects}
In this section we give some details of our implementation.  More
precisely, we explain some choices we made when implementing the
subroutines \verb/solve_over_Z_p/, \verb/eval_or_reduce/ and
\verb/check/ mentioned in the algorithm in
Section~\ref{sec:rational-reconstruction}.\footnote{The actual code
  can be found in the files {\tt modhpsol.spad.pamphlet} and {\tt
    mantepse.spad.pamphlet} in the \FriCAS\ distribution.}

\subsection{solving over $\mathbb Z_p$}

As already mentioned, our solver over $\mathbb Z_p$ closely
follows~\cite{BeckermannLabahn1994}.  It returns a matrix of
polynomials, every column constitutes one solution of the
Hermite-Pad\'e interpolation problem.  Solutions, and also the
residuals that occur within the algorithm, are packed in vectors
of machine size integers.  This is possible, since
we perform the same operations on 
each component polynomial of the solution.
However, it lowers control and memory management overhead.

Computations are performed in place -- otherwise memory management
would dominate the run time.  As basic operation we use \lq multiply
and add\rq, that is we compute $v_1 + cv_2$ for two vectors $v_1$ and
$v_2$ and a scalar $c$, and assign the result to $v_1$.  Compared to
separate addition and multiplication with a scalar, this approach halves the
cost of remainder computations needed for modular arithmetic and
also halves the loop overhead.

Currently the \lq multiply and add routine\rq\ is written in Spad (\FriCAS'
high level implementation language) and via Common
Lisp compiled to
machine code.\footnote{\FriCAS\ can compile code using a variety of
  Common Lisp implementations.  Currently, best performance is
  obtained with \sbcl, \url{http://www.sbcl.org}, thus in the
  following we assume this implementation.}
On 64-bit machines 32-bit times 32-bit multiplication
and 64-bit by 32-bit remainder are compiled directly to machine 
instructions.  On a 64-bit Core 2 this leads to about 20 clocks
per multiply and add step -- it seems that this is the same as
the cost of the machine instruction to compute a remainder.
On 32-bit machines our compilation scheme performs operations
involving 64-bit numbers (either as a result or as an argument)
via calls to bignum routines, which causes much higher
execution times.

In principle we could speed up the solver over $\mathbb Z_p$ replacing
remainder operations by multiplications.  Moreover, the \lq multiply
and add\rq\ routine is quite small so it would make sense to replace
it by an assembler routine.  We estimate that this would make
our inner loop run 5-10 times faster.  However, after
some initial work our solver over $\mathbb Z_p$ turned out to be the
fastest part of the package, so we concentrated
on removing bottlenecks in other parts.  Also, despite the
quadratic complexity and the non-optimal inner loop our solver
seems to compare favourably with other programs, like \GFUN,
\GuessKauers\ and \qGeneratingFunctions.

\subsection{computing modular images}
Initially we used a very naive evaluation scheme to obtain the
vector of truncated power series (or sequences) over $\mathbb Z_p$:
we did the computation in characteristic $0$ and then computed
remainders of division by $p$.  Measurements showed that this
approach used more than 98\% of the execution time.  Therefore, we
switched to a faster scheme.

For polynomials with integer coefficients
we use a specialised routine to reduce coefficients modulo $p$.
Also, for univariate polynomials with coefficients in $\mathbb Z_p$ we use
a specialised evaluation
routine.  
However, for multivariate polynomials with modular coefficients we still
use the naive method: we substitute for the variables in characteristic
$0$ and then use the routine just mentioned to reduce coefficients modulo $p$.
Of course this needlessly uses bignum arithmetic, but since multivariate
evaluation is typically followed by several univariate evaluations
the cost seem to be acceptable.  

Besides making the evaluation routines fast enough, it is important
to generate the vector $\Mat f$ of sequences or truncated power
series efficiently.  Initially we computed $\Mat f$ in characteristic
$0$ and then performed modular reduction.  It turned out that doing
the modular reduction only on the original sequence and computing the
derived sequences only over $Z_p$ is much faster.

Moreover, we remark that for truncated power series computing $\Mat
f$ involves computing many Cauchy products, which can be expensive.
Therefore we implemented a simple optimiser, that takes a vector of
monomials and detects common factors that can be cached.  This
reduced the number of Cauchy products that have to be computed
significantly.

With these improvements, the time needed for computing the modular
images typically is comparable to the time needed for solving over
$Z_p$.
\subsection{checking solutions}

For large problems, checking the
solutions may be dominant factor.  For sequences
operations are performed pointwise, but for truncated power
series we need polynomial multiplication which is much more
expensive -- the current polynomial routines of \FriCAS\ are quadratic.
Also, memory use may be a problem: we can solve the
Hermite-Pad\'e problem without actually computing the series
or sequences forming it (we only compute modular images), but
we need the actual system to check the solutions.  For example,
we can guess an equation for the generating function of Gessel
walks (see the article by Manuel Kauers, Christoph Koutschan, and
Doron Zeilberger about its holonomicity~\cite{MR2538821} and the
article by Alin Bostan and Manuel Kauers about its
algebraicity~\cite{BostanKauers2010})
using a few hundred megabytes of memory, but 
explicitly storing the Hermite-Pad\'e
problem needs several gigabytes and we run out of memory
on an 8~GB machine.  

We therefore introduced options, see Section~\ref{sec:options}, that allow 
skipping the checks or to use a Monte-Carlo check.  Although on small problems
the Monte-Carlo check may be slower than the deterministic check,
for the problem of Gessel walks it only takes about 30 sec.\ on a 2.4 GHz Core 2
and has moderate memory usage.

There is an additional problem: our initial data are fractions.
\FriCAS\ fraction arithmetic simplifies fractions after each
operation, which is quite costly (but also avoids intermediate
expression swell which would be even more costly).  In
practise in many cases denominator is $1$, and using fractions
causes almost no extra cost.  But some problems really make
use of fractional coefficients and for such systems checking
is more costly.  For linear recurrences we implemented a special
purpose checking routine which computes the common denominator
using smaller number of gcd operations and performs the rest
of operations in fraction-free way.  This routine also avoids
storing the vector $\Mat f$ of the Hermite-Pad\'e problem and 
thus avoids problems with excessive
memory usage.  Unfortunately, for other problems checking is
much more complicated so currently for them we use general
purpose checking routine which uses standard fraction
arithmetic.

We also considered using a modular method for checking.  Such a check
would avoid problems with memory use.  However, our a priori bounds
on the number of evaluation points needed are so
large, that the modular check is likely to be
slower then our current version.  Thus we only implemented a Monte-Carlo
version of a modular check (currently only for the case of
truncated power series).

\section{Normalisation}
The output of the algorithm proposed by Bernhard Beckermann and
George Labahn in 1994~\cite{BeckermannLabahn1994} is a so called
\Dfn{$\sigma$-basis} (also known as \Dfn{order basis}) for the
solution space of the rational interpolation problem.

In this section, we will call a polynomial vector a \Dfn{solution} of
the rational interpolation problem if it satisfies the order
condition~\eqref{eq:orderSequence} or \eqref{eq:orderSeries}.  The
degree constraints will be taken into account in a different manner,
namely through the notion of \Dfn{defect}:
\begin{dfn}
  Let $\Mat n = [n^{(1)},\dots,n^{(m)}]$ a vector of degree bounds
  and $\Mat p= [p^{(1)}(x),\dots,p^{(m)}(x)]$ be a vector of polynomials.

  Then the \Dfn{defect} of $\Mat p$ is the (possibly negative)
  difference between degree bound and degree of $\Mat p$. More
  precisely:
  $$\defect\Mat p=\min\left\{n^{(i)}-\deg p^{(i)}(x): 
    i\in\{1,\dots,m\}\right\}.
  $$

  For any $\delta\in\mathbb Z$, we denote by $\mathcal
  L^\sigma_\delta$ the space of solutions with defect strictly larger
  than $-\delta$.

  A \Dfn{$\sigma$-basis} $\{\Mat p_1,\dots,\Mat p_m\}$ is a set of
  solutions of the rational interpolation problem, such that every
  solution $\Mat q= [q^{(1)}(x),\dots,q^{(m)}(x)]$ can be written as
  a linear combination $\Mat q=\sum_{r=1}^m\alpha_r(x) \Mat p_r$,
  with $\defect\Mat q\leq\defect\Mat p_r-\deg\alpha_r(x)$, in a
  unique way.  Note that the elements of the $\sigma$-basis do not
  need to satisfy the degree bounds, i.e., their defect can be
  negative.
\end{dfn}

An alternative description of $\sigma$-bases is given as follows,
see~\cite[Equation~(6) and~(7)]{BeckermannLabahn1994}:
\begin{equation}
  \label{eq:spanningset}
  \mathcal L^\sigma_\delta = 
  \Span \{x^j\Mat p_r: 1\leq r\leq m, j<\defect \Mat p_r + \delta\},
\end{equation}
and
\begin{equation}
  \label{eq:dimension}
  \dim\mathcal L^\sigma_\delta = \sum_{r=1}^m \max(\defect\Mat p_r - d, 0).
\end{equation}
It follows that the set of defects is an invariant of the rational
interpolation problem.

We want to reconstruct a $\sigma$-basis over an integral domain $R$
(usually the integers or polynomials with integer coefficients) from
several $\sigma$-bases over quotient rings $R/I$, where $I$ is a
prime ideal in $R$.  For this purpose it is crucial to know when a
$\sigma$-basis over the field of fractions $Frac(R/I)$ of $R/I$ is an
image of the $\sigma$-basis over $R$.  The difficulty is that for a
given Hermite-Pad\'e problem, the $\sigma$-basis is not uniquely
determined.  Fortunately, over a field we may obtain uniqueness using
an appropriate normalisation.

In the following, we define such a normalisation and show the
existence and the uniqueness of a normalised $\sigma$-basis.  We then
prove a statement that can be applied to detect \lq bad
reductions\rq, i.e., to decide whether a given normalised
$\sigma$-basis over $Frac(R/I)$ is a modular image of a normalised
$\sigma$-basis over $Frac(R)$.

\begin{dfn}
  Let $\Mat p= [p^{(1)}(x),\dots,p^{(m)}(x)]$ be a polynomial vector,
  and $\Mat n = [n^{(1)},\dots,n^{(m)}]$ a vector of degree bounds.
  The \Dfn{critical index} of $\Mat p$ is the minimal index of the
  vector where the defect is attained:
  $$\critical\Mat p=\min\left\{i: 
    n^{(i)}-\deg p^{(i)}(x) = \defect\Mat p,
    i\in\{1,\dots,m\}\right\}.$$

  $\Mat p$ is \Dfn{normalised} if the leading coefficient of
  $p^{(i)}(x)$ is $1$, when $i$ is the critical index of $\Mat p$.

  $\Mat p$ is \Dfn{reduced} with respect to another polynomial vector
  $\Mat q= [q^{(1)}(x),\dots,q^{(m)}(x)]$, if at the critical index
  $i$ of $\Mat q$ we have $\deg p^{(i)}(x)<\deg q^{(i)}(x)$.  A
  sequence of polynomial vectors $\{\Mat q_1,\dots,\Mat q_m\}$ is
  \Dfn{reduced} if $\Mat q_r$ is reduced by $\Mat q_s$ for all $r\neq
  s$.  Note that this implies that the critical indices of the
  vectors must be all different.

  A sequence of polynomial vectors $\{\Mat q_1,\dots,\Mat q_m\}$ is
  \Dfn{sorted}, if for $r<s$
  \begin{align*}
    \defect\Mat q_r &> \defect\Mat q_s\quad\text{or}\\
    \defect\Mat q_r &= \defect\Mat q_s\quad\text{and} \critical\Mat
    q_r < \critical\Mat q_s.
  \end{align*}

  A sequence of polynomial vectors is \Dfn{normalised}, if it is
  sorted, reduced, and all its elements are normalised.
\end{dfn}

\begin{lem}
  Let $R$ be an integral domain.  If there is a normalised
  $\sigma$-basis over $R$ for 
  the solution space of a given rational interpolation problem,
  then it is uniquely determined.
\end{lem}
\begin{proof}
  Suppose we have two normalised bases $\Mat P=\{\Mat p_1,\dots,\Mat
  p_m\}$ and $\Mat Q=\{\Mat q_1,\dots,\Mat q_m\}$, and suppose that
  $\Mat p_r=\Mat q_r$ for $r<t$.  Assume without loss of generality
  that $d:=\defect\Mat q_t\geq\defect\Mat p_t$.  We will first show
  that the free module generated by the vectors in $\Mat P$ with
  defect $d$ coincides with the free module generated by the vectors
  in $\Mat Q$ with defect $d$.

  Consider any $\Mat q_r$ with $\defect\Mat q_r = d$, and expand it
  in the $\sigma$-basis $\Mat P$: $\Mat q_r = \sum_{s=1}^m
  \alpha_s(x) \Mat p_s$.  We want to show that $\deg\alpha_s(x)=0$,
  i.e., $\alpha_s(x)\in R$ for all $s$.

  Consider any $s$ with $\alpha_s\neq 0$.  Since $\Mat P$ is a
  $\sigma$-basis, we have $d=\defect\Mat q_r\leq\defect\Mat p_s$.
  Suppose that $\defect\Mat q_r<\defect\Mat p_s$.  If $s\geq t$ then
  $\defect\Mat p_s\leq\defect\Mat p_t$, since $\Mat P$ is sorted,
  contradicting the assumption that $d\geq\defect\Mat p_t$.  However,
  if $s<t$ we have $\Mat q_s=\Mat p_s$, and therefore that $\Mat q_r$
  is reduced with respect to $\Mat p_s$.  Thus, at the critical index
  $i$ of $\Mat p_s$ the degree of $q_r^{(i)}(x)$ is smaller than the
  degree of $p_s^{(i)}(x)$, which contradicts $\alpha_s\neq 0$.  So
  we must have $\defect\Mat q_r=\defect\Mat p_s$, and therefore
  $\deg\alpha_s(x)=0$.

  Similarly, expanding any $\Mat p_r$ with defect $d$ in the
  $\sigma$-basis $\Mat Q$ as $\Mat p_r=\sum_{s=1}^m \beta_s(x) \Mat
  q_s$ we obtain $\beta_s(x)\in R$: $\Mat Q$ is a $\sigma$-basis, so
  we have $d=\defect\Mat p_r\leq\defect\Mat q_s$.  Suppose that
  $\defect\Mat p_r<\defect\Mat q_s$.  Since $\Mat Q$ is sorted we
  must have $s<t$ and thus $\Mat p_s=\Mat q_s$.  We then conclude as
  above that $\deg\beta_s(x)=0$.

  This shows that the free modules over $R$ spanned by $\left\{\Mat
    p_{r_0},\dots,\Mat p_{r_1}\right\} := \{\Mat p_r: \defect\Mat p_r
  = d\}$ and by $\left\{\Mat q_{s_0},\dots,\Mat q_{s_1}\right\} :=
  \{\Mat q_s: \defect\Mat q_s = d\}$ coincide.  Note that $r_0=s_0$:
  because of $\defect\Mat p_{r_0}=d\geq\defect\Mat p_t$ we have
  $r_0\leq t$.  Thus, if we had $s_0<r_0$, then $s_0<t$ which implies
  $\Mat p_{s_0}=\Mat q_{s_0}$, and finally $\defect\Mat p_{s_0}=d$, a
  contradiction.  Suppose now $r_0<s_0$, then $r_0<t$, so $\Mat
  q_{r_0}=\Mat p_{r_0}$ and therefore $\defect\Mat q_{r_0}=d$, again
  a contradiction.  The number of vectors in both sequences must be
  the same, too, so $r_1=s_1$. 

  It remains to prove that the two polynomial sequences are in fact
  identical.
  If $r_0=r_1$, i.e., the module is one-dimensional, then the
  condition that $\Mat p_{r_0}$ and $\Mat q_{r_0}$ are normalised
  implies that they are identical.

  If $r_0<r_1$, we first show that the critical index $i$ of $\Mat
  p_{r_0}$ is the same as the critical index $j$ of $\Mat q_{r_0}$.
  Since $\Mat p_{r_0}$ is in the span of $\left\{\Mat
    q_{r_0},\dots,\Mat q_{r_1}\right\}$, there must also be a $\Mat
  q_r\in\left\{\Mat q_{r_0},\dots,\Mat q_{r_1}\right\}$ with $\deg
  q_r^{(i)}\geq n^{(i)}-d=\deg p_{r_0}^{(i)}$.  In fact, we have $\deg q_r^{(i)}=
  n^{(i)}-d$, since the defect of $\Mat q_r$ is $d$.  It follows that
  the critical index of $\Mat q_r$ is at most $i$, and, since
  $\left\{\Mat q_{r_0},\dots,\Mat q_{r_1}\right\}$ is sorted, $j\leq
  i$.  Interchanging the r\^oles of $\left\{\Mat p_{r_0},\dots,\Mat
    p_{r_1}\right\}$ and $\left\{\Mat q_{r_0},\dots,\Mat
    q_{r_1}\right\}$, we obtain $i\leq j$ and therefore $i=j$.

  Thus also the modules spanned by $\left\{\Mat
    p_{r_0+1},\dots,\Mat p_{r_1}\right\}$ and $\left\{\Mat
    q_{r_0+1},\dots,\Mat q_{r_1}\right\}$ coincide: we can obtain
  both from the module spanned by $\left\{\Mat p_{r_0},\dots,\Mat
    p_{r_1}\right\}$ by restricting to the set of polynomial vectors
  $\Mat v$ with $\deg v^{(i)}<n^{(i)}-d$.

  Iterating the argument of the previous two paragraphs, we find that
  the two modules generated by $\Mat p_{r_1}$ and $\Mat q_{r_1}$
  coincide, and since these vectors are normalised, they must be
  identical.  In particular, their critical indices coincide.
  Because the sequences $\left\{\Mat p_{r_0},\dots,\Mat
    p_{r_1}\right\}$ and $\left\{\Mat q_{r_0},\dots,\Mat
    q_{r_1}\right\}$ are reduced, we can reuse the argument of the
  previous paragraph with $i$ being this critical index, to obtain
  that also the modules spanned by $\left\{\Mat p_{r_0},\dots,\Mat
    p_{r_1-1}\right\}$ and $\left\{\Mat q_{r_0},\dots,\Mat
    q_{r_1-1}\right\}$ coincide.
  Now induction shows that $\Mat p_r=\Mat q_r$ for
  $r\in\{r_0,\dots,r_1\}$, which concludes the proof of uniqueness.
\end{proof}

\begin{lem}
  The solution space of
  every rational interpolation problem has a normalised
  $\sigma$-basis over a field.
\end{lem}
\begin{proof}
  The existence of $\sigma$-bases over a field is meanwhile well
  known, for example Bernhard Beckermann and George Labahn prove that
  their algorithm produces a $\sigma$-basis
  in~\cite{BeckermannLabahn1994}.  Let $\Mat P=\{\Mat p_1,\dots,\Mat
  p_m\}$ be a $\sigma$-basis.  We will show that for any $d$ we can
  replace $\{\Mat p_r: \defect\Mat p_r \geq d\}$ by an equivalent
  normalised sequence.  We proceed by induction: assume that $S_0 =
  \{\Mat p_1,\dots\Mat p_k\}$ is a normalised sequence of vectors
  with $\defect\Mat p_r > d$, and let $S_1=\{\Mat p_{k+1},\dots,\Mat
  p_l\}$ is a sequence of vectors with defect $d$, which we will
  successively add to $S_0$.

  First of all we remark that replacing any vector $\Mat q$ of a
  $\sigma$-basis $\Mat P$ by $\Mat{\tilde q}=\Mat q+\alpha\Mat p$,
  where $\Mat p\in\Mat P$ and $\defect\Mat p\geq\defect\Mat
  q+\deg\alpha$, again yields a $\sigma$-basis.  Moreover, since by
  Equation~\eqref{eq:dimension} the set of defects is an
  invariant of the solution space, we have $\defect\Mat{\tilde
    q}=\defect\Mat q$.

  Thus, by subtracting appropriate polynomial multiples of previous
  elements we may assume that each $\Mat p_t\in S_1$ is reduced with
  respect to $\Mat p_r \in S_0$.  Namely,
  let $\Mat q \in S_1$,
  let $Q_0$ be set of critical indices of elements of $S_0$ and
  let $d_i=\defect\Mat p_r$ with $r$ such that $i$ is critical
  index of $\Mat p_r$.  Consider $c_i = n^{(i)}-\deg q^{(i)}$.
  If $c_i > d_i$ for all $i \in Q_0$, then $\Mat q$ is reduced with
  respect to all $\Mat p_r \in S_0$, $r < t_0$.  Otherwise consider
  $Q_1 = \{i : c_i \leq d_i\}$ and let $Q_2$ be set of $i \in Q_1$
  such that $c_i$ is minimal.
  Select the first
  element $\Mat p \in S_0$, with critical index $i_0 \in Q_2$.
  Then choose $\alpha$, such that
  $\deg (q^{(i_0)}-\alpha p^{(i_0)}) < \deg q^{(i_0)}$
  and replace $\Mat q$ in $S_1$ by 
  $\Mat{\tilde q} = \Mat q - \alpha \Mat p$.  Let 
  $b_i = n^{(i)} - \deg{\tilde q}^{(i)}$.  Note
  that $\deg \alpha = d_{i_0} - c_{i_0}$ and (by the definition of defect)
  $n^{(i)} - \deg p^{(i)} \geq d_{i_0}$ for all $i$, so
  $\deg(\alpha p^{(i)}) \leq n^{(i)} - c_{i_0}$.  Moreover, if $i$ is the
  critical index of $\Mat p_r$ which is smaller than $\Mat p$, then the
  defect of $\Mat p_r$ is bigger or equal to the defect of $\Mat p$, and,
  since $\Mat p$ is reduced with respect to $\Mat p_r$, we have
  $n^{(i)} - \deg p^{(i)} > d_{i_0}$ and
  $\deg(\alpha p^{(i)}) < n^{(i)} - c_{i_0}$.  Next
  $\deg{\tilde q})^{(i)} \leq \max\big(\deg q^{(i)}, \deg(\alpha p^{(i)}\big)$.
  Consequently, we have for all $i$
  \begin{align*}
    b_i &= n^{(i)} - \deg{\tilde q}^{(i)} \\
        &\geq n^{(i)} - \max\big(\deg q^{(i)}, \deg(\alpha p^{(i)})\big)\\
        &= \min\big(n^{(i)} - \deg q^{(i)}, n^{(i)} - \deg(\alpha p^{(i)})\big)\\
        &= \min\big(c_i, n^{(i)} - \deg(\alpha p^{(i)})\big)\\
        &\geq \min(c_i, c_0),
  \end{align*}
  and
  for $i$ corresponding to $\Mat p_r$ less or equal to $\Mat q$ we have
  $b_i \geq \min(c_i, c_0 + 1)$.  This means that after a finite
  number of reductions passing trough $\Mat q_0 = \Mat q$, 
  $\Mat q_1 = \Mat{\tilde q}$, $\Mat q_2$, etc. we will increase 
  $\min\{n^{(i)} -  \deg q_j^{(i)}): i \in Q_1 \}$.
  Continuing the reduction process we will increase 
  $\size{\{i \in Q_0 : n^{(i)} -  \deg q_j^{(i)} > d_i \}}$ and
  eventually $\Mat q_j$ will be reduced with respect to all elements
  of $S_0$.

  Note that every element $\Mat p_r\in S_0$ is automatically reduced
  with respect to every element $\Mat q\in S_1$, because the defect
  of $\Mat p_r$ is strictly smaller than the defect of $\Mat q$: let
  $i$ be the critical index of $\Mat p_r\in S_0$, then $n^{(i)}-\deg
  q^{(i)}\geq \defect\Mat q > \defect\Mat p_r = n^{(i)}-\deg
  p_r^{(i)}$.

  It remains to show how to ensure that $\Mat p_r\in S_1$ is reduced
  with respect to $\Mat p_s\in S_1$.  This can be
  achieved by applying a simpler version of process described above:
  in a first step
  reduce $\Mat p_{k+2}, \dots, \Mat p_l$ with respect to
  $\Mat p_{k+1}$, then $\Mat p_{k+3}, \dots, \Mat p_l$ with respect to the
  sequence $\{\Mat p_{k+1}, \Mat p_{k+2}\}$ and so on.  Note that,
  since the defects of these vectors are all equal, the degree of
  $\alpha$ will always be zero.  This implies that the vector
  replacing $\Mat p_s$ will remain reduced with respect to $\Mat
  p_r$, $r<s$.  

  After this step, all $\Mat p_s\in S_1$ are
  reduced with respect to $\Mat p_r\in S_1$ with $r < s$.
  Similarly, but going backwards, we ensure that
  $\Mat p_s\in S_1$ are
  reduced with respect to $\Mat p_r\in S_1$ with $r > s$.
\end{proof}

To continue we need a lemma about linear systems.
\begin{lem}\label{lem:bad-reduction}
  Let $n_1 = \dim \ker(A)_{Frac(R)}$ and $n_2 =
  \dim\ker(A)_{Frac(R/I)}$.  Then $n_1\leq n_2$.  Moreover, if $n_1 =
  n_2$, then $\ker(A)_{Frac(R/I)} = \pi\big(\ker(A)_{R_I}\big)$,
  where $\pi$ is the quotient map and $R_I$ is the localisation of
  $R$ at $I$, i.e., the ring of fractions $\frac{r}{s}$ with $r\in R$
  and $s\in R\setminus I$.
\end{lem}
\begin{proof}
For vector spaces we have 
$\dim\ker(A) = \dim\dom(A) - \rank(A)$, so $n_1 \leq n_2$
is equivalent to $\rank(A)_{Frac(R)} \geq \rank(A)_{Frac(R/I)}$.
This inequality follows easily by considering
minors of $A$.  So it remains to prove that $n_1 = n_2$ implies 
$\ker(A)_{Frac(R/I)} = \pi(\ker(A)_{R_I}$.  Let
$m = \rank(A)_{Frac(R/I)}$.  It is enough
to prove this equality for surjective $A$ over $Frac(R)$.
Namely permuting rows of $A$ we can write $A$ in block form:
$$
\left(
    \begin{array}{c}
     A_1 \\
     A_2 \\
     \end{array}
 \right)
$$
such that $A_1$ has $m$ rows and $\rank(A_1)_{Frac(R/I)} = m$.  Since
$n_1 = n_2$ this means
that also $\rank(A_1)_{Frac(R/I)} = \rank(A_1)_{Frac(R)}$.
Next, $\ker(A)_{Frac(R)} \subset \ker(A_1)_{Frac(R)}$, and
$\dim\ker(A)_{Frac(R)} = \dim\ker(A_1)_{Frac(R)}$, so
$\ker(A)_{Frac(R)} = \ker(A_1)_{Frac(R)}$.  Consequently,
$A_2 = 0$ on $\ker(A_1)_{Frac(R)}$, so also
$\ker(A)_{R_I} = \ker(A_1)_{R_I}$.  Similarly,
$\ker(A)_{Frac(R/I)} = \ker(A_1)_{Frac(R/I)}$.  So,
indeed it is enough to prove
$\pi\big(\ker(A_1)_{R_I}\big) = \ker(A_1)_{Frac(R/I)}$
and replacing $A$ by $A_1$ we may assume that $A$ is
surjective.
By permuting columns of $A$ we may assume that $A$
has block form:
$$
\left(
    \begin{array}{cc}
     A_1 & A_2 \\
     \end{array}
 \right)
$$
where $A_1$ is invertible over $Frac(R/I)$.  This
means that determinant of $A_1 \notin I$, so $A_1$ is
invertible over $R_I$.  Multiplying from the left
by $A_1^{-1}$ we may assume that $A_1$ is identity
matrix.  But then it is easy to compute the kernel over
$R_I$: it consists of the vectors of the form $(-A_2v, v)$,
where $v$ is an arbitrary vector in the domain of $A_2$.
This implies that
$\dim\pi(\ker(A)_{R_I} = \dim\ker(A)_{Frac(R/I)}$,
which gives the claim.
\end{proof}

\begin{cor}\label{cor:bad-reduction} 
  Let $\Mat P$ be a $\sigma$-basis over $Frac(R)$ and let $\Mat Q$ be 
  a $\sigma$-basis over $Frac(R/I)$.  Assume that $\defect\Mat p_r \geq
  \defect\Mat p_{r+1}$ and $\defect\Mat q_r\geq\defect\Mat q_{r+1}$ for all $r$.  Then
  $\defect\Mat p_r \leq \defect\Mat q_r$ for all $r$.  If both 
  $\Mat P$ and $\Mat Q$ are normalised and the defects and critical indices of $\Mat p_r$ and 
  $\Mat q_r$ are the same for all $r$, then $\Mat P$ is defined over $R_I$ and $\Mat Q$ is an 
  image of $\Mat P$ via the quotient map $\pi$.
\end{cor}
\begin{proof}
  According to Equation~\eqref{eq:dimension}, the dimensions of the
  solution spaces $\mathcal L^\sigma_\delta$ of the rational interpolation problem over
  $Frac(R)$ and $Frac(R/I)$ are given by $\sum_{r=1}^m
  \max(\defect\Mat p_r-\delta, 0)$, and $\sum_{r=1}^m
  \max(\defect\Mat q_r-\delta, 0)$ respectively.  By
  Lemma~\ref{lem:bad-reduction}, and choosing $\delta\in\mathbb Z$ small enough,
  we obtain that $\defect\Mat p_r\leq\defect\Mat q_r$.

  When all defects are equal, the dimensions of the solution spaces
  coincide, so by the second part of the lemma, the solutions over
  $Frac(R/I)$ are images of solutions over $R_I$.
  In particular, for each $r$ the solution $\Mat q_r$ is an image of
  a solution $\Mat h_r$ over $R_I$.  By
  Equations~\eqref{eq:spanningset} and~\eqref{eq:dimension}, 
  the sequence $\Mat H = \{\Mat
  h_1,\dots,\Mat h_m\}$ is a $\sigma$-basis over $Frac(R)$.  It is
  easy to see that normalising $\Mat H$ gives a $\sigma$-basis which
  is equal to $\Mat H\mod I$ (here we need the assumption on the
  critical indices).  Since there is a unique normalised
  $\sigma$-basis, $\Mat H\mod I$ and $\Mat Q$ coincide.
\end{proof}

\begin{cor} The previous corollary remains valid if
we only assume that $I$ is an intersection of prime ideals,
and that $\sigma$-basis computation and normalisation
worked over $Frac(R/I)$ without encountering division by
non-invertible element.
\end{cor}
\begin{proof}
If $I=I_1\cap\dots\cap I_n$ then $Frac(R/I)$ is isomorphic
to a subring of product $\Pi_{i=1}^nFrac(R/I_i)$ and the claim easily follows.
\end{proof}

Let us stress that the corollary above means that either all modular
images are \lq bad reductions\rq\ (that is the dimension of the space
of modular solutions is bigger than the dimension of the original
space, or not all critical indices are equal) which is highly
improbable, or, by normalising and rejecting solutions with larger
defect or different critical indices we obtain a consistent
normalisation.

\section{Performance}
\label{sec:performance}

To test the performance of the package, we ran a few examples with
our package, \GFUN\ (version~3.5 on Maple~11), and \GuessKauers\
(version~0.32 on Mathematica~7.0).

Timings are in seconds, best of three runs.  \GuessKauers\ was run on
a {\tt Intel Core 2 E8400 @ 3 GHz} with 6MB cache and 1.8GB RAM but
running a 32-bit operating system, the other two on a {\tt Intel
  Pentium 4 @ 3 GHz}, 2MB cache, 1 GB RAM.

Both \Guess\ and \GFUN\ tried all configurations of order and degree,
only \GuessKauers\ was run with specified order and degree of the
recurrence.  Since both \GFUN\ and \GuessKauers\ look for homogeneous
recurrences by default, we invoked \spad{guessPRec} with
\spad{homogeneous==true}.  We believe that neither \GFUN\ not
\GuessKauers\ check the recurrence found, thus \spad{guessPRec} was
invoked with \spad{check=='skip}.

On the one hand, we recovered randomly generated homogeneous
polynomial recurrences over $\mathbb Q$ from data, see
Table~\ref{tab:randomrec}.  On the other hand, we computed
homogeneous polynomial recurrences over $\mathbb Q[t]$ for the first
few integer powers of the Hermite polynomials, see
Table~\ref{tab:Hermiterec}.  (This second test was
only run against \GFUN.  For comparison, we also indicate order and
degree of the recurrence discovered.)

Readers should be cautious interpreting the data.  Theoretically, the
algorithm used by \GFUN\ has lower complexity for large degrees,
while \GuessKauers\ seems best adapted to very low degrees.  However,
as we explained performance depends very much on implementation
detais and we lack sufficient information about the other packages to
make a more general and precise statement.

\begin{table}
  \centering
  $$
  \begin{array}[h]{r|r|r|r|r|r|r}
    \text{order}\backslash\text{degree}
    &  10 &   20 &  30 &   40 &    50   \\\hline\hline
  5 &  0.0&   0.1&  0.3&   0.6&    1.1  \\
    &  0.3&   1.1&  1.7&   4.3&    5.7  \\
    &  0.1&   0.5&  1.8&   5.0&   11.3  \\\hline
 10 &  0.1&   0.6&  1.6&   3.6&    7.4  \\   
    &  1.7&   5.6&  8.8&  20.9&   27.8  \\
    &  0.3&   2.4& 10.8&  29.2&   65.0  \\\hline
 15 &  0.4&   2.0& 5.4 &  12.5&   24.5  \\   
    &  5.6&  16.0& 42.6&  59.4&   77.3  \\
    &  0.9&   7.5& 33.3&  87.1&  201.6  \\\hline
 20 &  1.0&   4.8& 13.8&  31.6&  115.6  \\  
    & 19.2&  39.4& 99.0& 137.1&  179.0  \\
    &  1.9&  15.3& 69.4& 196.5&  447.5  \\\hline
 25 &  2.2&  10.2&     &      &         \\
    & 40.3&  85.3&     &      &         \\
    &  3.3&  30.8&     &      &         \\\hline
 30 &  4.2&  18.7&     &      &         \\
    & 75.8& 162.8&     &      &         \\
    &  5.2&  50.1&     &      &         \\\hline
 35 & 7.3 &  32.6&     &      &         \\
    &132.3& 278.6&     &      &         \\
    & 7.7 &  75.8&     &      &         \\\hline
 40 & 11.6&  51.8&     &      &         \\
    &221.0& 604.9&     &      &         \\
    & 11.0& 120.0&     &      &         \\\hline
 45 & 17.8&  78.1&     &      &         \\
    &353.6& 906.5&     &      &         \\
    & 15.1& 164.9&     &      &         \\\hline
 50 & 26.1& 116.5&     &      &         \\
    &536.1& 309.4&     &      &         \\
    & 20.0& 219.7&     &      &         \\\hline
 55 & 37.0& 163.1&     &      &         \\
    &787.3& 838.7&     &      &         \\
    & 26.0& 285.6&     &      &         \\\hline
 60 & 52.3& 222.2&     &      &         \\
   &1094.6&2516.6&     &      &         \\
    & 38.6& 362.3&     &      &         \\\hline
 65 & 70.3&      &     &      &         \\
   &1509.0&      &     &      &         \\
    & 48.0&      &     &      &         \\\hline
 70 & 92.5&      &     &      &         \\
   &1927.2&      &     &      &         \\
    & 58.5&      &     &      &         \\\hline
    \end{array}
    $$
    \caption{Guessing random homogeneous recurrences with polynomial
      coefficients over $\mathbb Q$.  The first line is \Guess, the
      second \GFUN, the third \GuessKauers.}
    \label{tab:randomrec}
\end{table}

\begin{table}
  \centering
  $$
\begin{array}[h]{r|r|r|r|r|r|r}
       & H(.,t)^1  & H(.,t)^2& H(.,t)^3 & H(.,t)^4 &   H(.,t)^5 & H(.,t)^6\\\hline\hline
\text{order}&  3   &  4 &   5 &   6 &     7 & 8\\
\text{degree}&  1   &  3 &   7 &  13 &    22 & 34\\\hline
\Guess &  0.0 & 0.0&  0.0&  0.4&   5.1 & 46.2\\
\GFUN  &  0.0 & 0.1&  2.5& 20.2& 238.3 & \text{fail}\\
\end{array}
$$
    \caption{Guessing random homogeneous recurrences for powers of
      the Hermite polynomials. (\GFUN\ ran out of memory computing
      the last entry.)}
    \label{tab:Hermiterec}
\end{table}

\section{Further work}

To conclude, we would like to point out possible future directions:
\begin{itemize}
\item It would be very important to generalise to the
  multidimensional case, as already implemented by Manuel Kauers in
  his package.  Of course, we can employ \lq diagonal guessing\rq,
  see~\cite{Zeilberger2007}.  I.e., we could first guess formulas for
  each row, and then guess formulas for the coefficients of these.
  However, this approach is rather slow and, more importantly,
  depends on the availability of many terms.

\item The performance of \spad{guessExpRat} and \spad{guessBinRat} is
  very disappointing, making the two procedures nearly useless.
  Moreover, these two are but a toy example for real world
  applications, where one would like to guess formulas like
  \begin{equation*}
    \det\left(\binom{3(i+j)+1}{i+j}\right) =
    \prod_{i=1}^n\frac{(6i+4)!(2i+1)!}{2(4i+2)!(4i+3)!}
    \sum_{i=0}^n\frac{n!(4n+3)!!(3n+i+2)!}{(3n+2)!i!(n-i)!(4n+2i+3)!!},
  \end{equation*}
  as found by {\"O}mer E{\u{g}}ecio{\u{g}}lu, Timothy Redmond and
  Charles Ryavec~\cite{MR2557881}.

\item Maybe there are other interesting operators besides $\Delta_n$ and $Q_n$
  that could be applied recursively to the sequence. Furthermore, there is a
  list of transformations used in \emph{The online encyclopedia of integer
    sequences}, it might be rewarding to check which of those extend the class
  of functions already covered significantly.
\end{itemize}

\providecommand{\cocoa} {\mbox{\rm C\kern-.13em o\kern-.07em C\kern-.13em
  o\kern-.15em A}}
\providecommand{\bysame}{\leavevmode\hbox to3em{\hrulefill}\thinspace}
\providecommand{\MR}{\relax\ifhmode\unskip\space\fi MR }
\providecommand{\MRhref}[2]{%
  \href{http://www.ams.org/mathscinet-getitem?mr=#1}{#2}
}
\providecommand{\href}[2]{#2}

\end{document}